\input amstex
\documentstyle{amsppt}
\magnification=1200
\NoBlackBoxes
\NoRunningHeads
\document

\topmatter
\title Nonabelian Integrable Systems, Quasideterminants,
and Marchenko Lemma
\endtitle

\author Pavel Etingof, Israel Gelfand, and Vladimir Retakh
\endauthor

\address
\newline
I.~G. :  Department of Mathematics, Rutgers University, New Brunswick,
NJ 08903
\newline
P.~E. : Department of Mathematics, Harvard University, Cambridge, MA
02138
\newline
V.~R. : Department of mathematics, University of Arkansas, Fayetteville, 
AR 72701
\endaddress

\email
\newline
P.~E. : etingof\@ math.harvard.edu
\newline
I.~G. : igelfand\@ math.rutgers.edu
\newline
V.~R. : vretakh\@ comp.uark.edu
\endemail

\endtopmatter

\centerline{\bf Abstract}

We find explicit (multisoliton) solutions for nonabelian integrable 
systems such as periodic Toda
field equations, Langmuir equations, and Schr\"odinger equations for
functions with values in any associative algebra. The solution
for nonabelian Toda field equations for root systems of types
$A, B, C$ was expressed by the authors in [EGR] using quasideterminants
introduced and studied in [GR1-GR4]. To find multisoliton solutions
of periodic Toda
equations and other nonabelian systems we use a combination of these
ideas with important lemmas which are due to Marchenko [M].

\centerline{\bf Introduction}
 
Let $R$ be an algebra. A map $\partial :R\rightarrow R$ is
called a derivation of $R$ if 
$\partial (ab)=(\partial a)b+a\partial b$.

Suppose that $R$ is an algebra with unit,
and $\partial _1$, $\partial _2$ two derivations of $R$
commuting with each other.

\proclaim {Definition 1} One says that invertible elements
$g _i, i=0,1,\dots , n-1$ satisfy the $n$-periodic Toda
equations if
$$
\partial _1((\partial _2g _k)g^{-1}_k)=
g _kg ^{-1}_{k-1} - g _{k+1}g^{-1} _k,
\  k\in {\Bbb Z}/n{\Bbb Z} .\tag 1  
$$
\endproclaim
In this paper we construct solutions of nonabelian periodic Toda systems
and other integrable systems using a combination of the theory of 
quasideterminants developed 
in [GR1, GR2, GR3, GR4] and Marchenko approach [M] to integrable
systems of differential equations on functions with values in 
operator algebras.

 Let $Q$
be a not necessarily commutative algebra with derivations 
$\partial _1,\dots ,
\partial _m$ and automorphisms $\alpha _1,\dots ,\alpha _s$.
Marchenko approach says that when an element $\Gamma \in Q$
satisfies a certain system of linear differential equations 
connecting the $\partial $'s, $\alpha $'s and $\Gamma $, then
its (noncommutative) logarithmic derivative 
$(\partial _1\Gamma)\Gamma ^{-1}$ satisfies a system of
nonlinear equations. 

Because Marchenko's book is unfortunately unavailable in English,
we will reproduce parts of its material for reader's convenience.

Note that nonabelian Toda equations for the root system $A_{n-1}$
were introduced by Polaykov (see [Kr]). Nonabelian Toda lattice for
functions of one variable appeared in [PC], [P]. Krichever [Kr] 
constructed
algebraic-geometric solutions for the periodic nonabelian Toda
lattices for $N\times N$-matrix-valued functions. 

\remark{Remark}
Our construction of multisoliton solutions is similar to the ideas
of higher rank solitons proposed by
Krichever and Novikov at the end of seventies (see e.g. \cite{KN}).
\endremark
    
\subhead 1.Quasideterminants and Wronski Matrices\endsubhead
Quasideterminants were introduced in \cite{GR1}, as follows. 
Let $X$ be an $m\times m$-matrix over $A$. 
For any $1\le i,j\le m$, let 
$r_i(X)$, $c_j(X)$ be the i-th row and the j-th column of $X$. 
Let $X^{ij}$ be the submatrix of $X$ obtained by removing 
the i-th row and the j-th column from $X$. For a row 
vector $r$ let $r^{(j)}$ be $r$ without the j-th entry.
For a column vector $c$ let $c^{(i)}$ be $c$ without the i-th entry. 
Assume that $X^{ij}$ is invertible. Then {\it the quasideterminant} 
$|X|_{ij}\in A$ is defined by the formula
$$
|X|_{ij}=x_{ij}-r_i(X)^{(j)}(X^{ij})^{-1}c_j(X)^{(i)},
$$
where $x_{ij}$ is the $ij$-th entry of $X$.

Let $Q$ be an algebra with unit, and 
$\partial :Q\rightarrow Q$ its derivation. For $f\in Q$ denote 
$\partial ^kf$ by $f^{(k)}, k=1,2,\dots $ and set
$f^{(0)}=f$.

For elements $f_1,\dots ,f_N$ denote by $W=W(f_1,\dots , f_N)$ their Wronski
matrix, 
$$
W=\left(\matrix f_1&...&f_N\\
\partial f_1&...&\partial f_N\\
..&...&..\\
\partial ^{N-1}f_1&...&\partial ^{N-1}f_N\endmatrix\right) .
$$
Set 
$\partial W=(\partial f^{(k-1)}_j)=(f^{(k)}_j0$,
$W^i=(f^{(m)}_j)$, $m=0,\dots ,\hat i,\dots , N$,
$j=1,\dots ,N$. Suppose that $W$ is invertible and
consider matrix $X=(\partial W)W^{-1}=(\eta _{pq})$.
Let $\delta _{pq}$ be the Kronecker symbol.
\proclaim {Proposition 1} For $p=1,\dots ,N-1$
$${\text i)}\ \eta _{pq}=\delta _{p+1,q},$$ 
and for $q=1,2,\dots ,N$
$$
{\text ii)}\ \eta _{Nq}=|W^q|_{NN}|W|^{-1}_{qN}.
$$
\endproclaim
{\bf Proof}. By our construction $\partial W=XW$.
It is easy to see that this identity is satisfied
if and only if i) is valid for $p=1,\dots ,N-1$,
and for $q=1,\dots ,N$
$$
f^{(N)}_q=\eta _{N1}f^{(0)}_q+\dots 
+\eta _{NN}f^{(N-1)}_q. \tag 2 
$$
Formulas ii) follow from Cramer rules [GR1, GR2, GR3, GR4]
applied to the right system of linear equations (2).
$\square $
\remark{Remark} Consider the $(N+1)\times N$-matrix
$$
\tilde W=\left(\matrix f_1&...&f_N\\
\partial f_1&...&\partial f_N\\
..&...&..\\
\partial ^{N}f_1&...&\partial ^{N}f_N\endmatrix\right) .
$$
Then $\eta _{qN}$ is the right quasi-Pl\"ucker
coordinate $r^{1\dots \hat i\dots N-1}_{Nq}(\tilde W)$
defined in [GR4].
\endremark
Square matrices $(\eta _{pq})$ satisfying Property i) 
in Proposition 1 are often called {\it Frobenius cells}.

\subhead 2. Marchenko Lemma\endsubhead    
Let $R$ be a ring with unit, $\partial _i : R\rightarrow R$,
$i=1,2$ its derivations, and $\alpha :R\rightarrow R$ an
isomorphism, such that $\partial _i\alpha =\alpha \partial _i$
for $i=1,2$.

Let also $\Gamma , A\in R$ be invertible elements such that
$\partial _iA=0, i=1,2$, $\partial _1(\partial _2\Gamma )=\Gamma$,
and $\partial _2\Gamma = \alpha (\Gamma)A$.

Set $\gamma = (\partial _2 \Gamma)
\Gamma ^{-1}$. The following lemma was proved in [M].  
\proclaim {Marchenko Lemma} 
$$
\partial _1(\partial _2 \gamma)\gamma^{-1}=
\gamma\alpha ^{-1} (\gamma^{-1})-\alpha (\gamma)\gamma^{-1}. \tag 3
$$
\endproclaim

{\bf Proof}. We will use the identity
$\delta (a^{-1})=-a^{-1}(\delta a)a^{-1}$ for an
invertible element $a\in R$ and a derivation $\delta :
R\rightarrow R$.  

One could write $\gamma $ as 
$\gamma =\alpha (\Gamma) A \Gamma ^{-1}$. Then

$$
\partial _2\gamma = (\partial _2\alpha (\Gamma ))A\Gamma ^{-1}
-\alpha (\Gamma)A\Gamma^{-1} (\partial _2\Gamma)\Gamma^{-1}=
$$
$$
=\alpha ^2(\Gamma)\alpha (A) A\Gamma ^{-1} - 
\alpha (\Gamma) A\Gamma^{-1}\alpha (\Gamma)A\Gamma^{-1}.
$$
Then
$$
(\partial _2\gamma)\gamma ^{-1}=
\alpha ^2(\Gamma)\alpha (A)\alpha(\Gamma)^{-1}
-\alpha(\Gamma)A\Gamma ^{-1}. \tag 4 
$$
Note that ($\partial _1\alpha (\Gamma ))A
=\partial_1 (\alpha (\Gamma )A)
=\partial _1(\partial _2\Gamma)=\Gamma.$ So,
$\partial _1(\alpha (\Gamma ))=\Gamma A^{-1}$, and
$\partial _1(\Gamma )=\alpha ^{-1} (\Gamma)\alpha ^{-1}(A)^{-1}$.

>From (4) it follows that
$$
\partial _1(\partial _2(\gamma )\gamma ^{-1})
=\alpha (\Gamma )\alpha (A)^{-1}\alpha (A)\alpha (\Gamma )^{-1}
-\alpha ^2(\Gamma)\alpha(A)\alpha (\Gamma )^{-1}
\Gamma A^{-1}\alpha (\Gamma )^{-1}-
$$
$$
-\Gamma A^{-1}A\Gamma ^{-1}
+\alpha(\Gamma )A\Gamma ^{-1}\alpha ^{-1}(\Gamma )\alpha ^{-1}(A)^{-1}
\Gamma ^{-1}=
$$
$$
=-\alpha ^2(\Gamma)\alpha(A)\alpha (\Gamma )^{-1}
\Gamma A^{-1}\alpha (\Gamma )^{-1}
+\alpha(\Gamma )A\Gamma ^{-1}\alpha ^{-1}(\Gamma )\alpha ^{-1} 
(A)^{-1}\Gamma ^{-1}.
$$

On the other hand, 
$$
\gamma\alpha ^{-1} (\gamma ^{-1}) - \alpha (\gamma )\gamma ^{-1}
=\alpha(\Gamma )A\Gamma ^{-1}\alpha ^{-1}(\Gamma )\alpha ^{-1} (A)^{-1}
\Gamma ^{-1}-
$$
$$
-\alpha ^2(\Gamma)\alpha(A)\alpha (\Gamma )^{-1}
\Gamma A^{-1}\alpha (\Gamma )^{-1}.
$$
$\square$
\subhead 3. Specialization of Marchenko Lemma\endsubhead
We will apply Marchenko Lemma to the following data.
Let $Q$ be an algebra, and 
let $\partial _i:Q\rightarrow Q, i=1,2$ be derivations
of $Q$ commuting with each other. Set 
$R=Q^n=\{(w_0,\dots ,w_{n-1})\}$ be an algebra with
coordinate-wise multiplication and coordinate-wise derivations
induced by $\partial _i, i=1,2$. These derivations will be 
denoted by the same letters.

To construct an isomorphism $\alpha :R\rightarrow R$,
let us suppose that indices $0,1,\dots ,n-1 \in {\Bbb Z}/n{\Bbb Z}$
 and set $\alpha ((w_k))=(w_{k+1})$. Then $\alpha ^{-1} ((w_k))=
(w_{k-1})$.

Let $A_k\in Q$, $k=0,1,\dots ,n-1$ be elements 
such that for $k=0,1,\dots , n-1$
$$\partial _iA_k=0,\ i=1,2$$
and $(w_k)$ be such that
$$
\partial _2w_k=w_{k+1}A_k, \ \partial _1w_k=w_{k-1}A_{k-1}^{-1}.
$$
Set $\gamma _k=(\partial _2w_k)w_k^{-1}$, $k=0,1,\dots ,n-1$.
Then Marchenko Lemma implies
\proclaim {Proposition 2}
$$
\partial _1((\partial _2\gamma _k)\gamma^{-1}_k)=
\gamma _k\gamma ^{-1}_{k-1} - \gamma _{k+1}\gamma^{-1} _k,
\  k=0,1,\dots , n-1. \tag 5 
$$
\endproclaim
\subhead 4. Frobenius cells\endsubhead
We will apply the results
of subsection 3 to Frobenius cells.

Consider matrices $K=(\mu _{pq})$, $L=(\nu _{pq})$, and $Y=(y_{pq})$,
$p,q=1,\dots , N$ over an algebra with unit.
Suppose that matrices $K$ and $L$ are Frobenius cells, i.e.
that $\mu _{pq}=\nu _{pq}=\delta _{p+1q}$
for $p=1,\dots ,N-1, q=1,\dots , N$ ($\delta $ is the Kronecker symbol).

\proclaim {Lemma 3} Suppose that $K=YL$ and that $\nu _{N1}$ is
invertible. Then
$y_{pq}=\delta _{pq}$ for $p=1,\dots , N-1, q=1,\dots , N$ and
$$
y_{Nq}=\mu _{N,q+1}-\mu _{N1}\nu ^{-1}_{N1}\nu _{N,q+1},
\ q=1,\dots , N-1,
$$
$$
y_{NN}=\mu _{N1}\nu ^{-1}_{N1}.
$$
\endproclaim
{\bf Proof}. Straightforward. $\square$

Suppose that under conditions of Lemma 3 
$\gamma _k=(g_{k,pq}), p,q=1,\dots , N,$ is a Frobenius
cell for $k=0,1,\dots, n-1$. Then Lemma 3 immediately implies
\proclaim{Corollary 4} Elements $g_k=g_{k, N1}$,
$k=0,1,\dots ,n-1$ satisfy the Toda periodic system (1).
\endproclaim
\subhead 5. Solutions of the Toda system\endsubhead 
Let $Q$ be an algebra with unit, and
let $\partial _1, \partial _2:
Q\rightarrow Q$ be two commuting derivations.

Let $f_{ij}, a_{ij}$
$i\in {\Bbb Z}/n{\Bbb Z}, j=1,\dots , N$ be elements in $Q$ such that
for all $i,j$ $a_{ij}$ are invertible and
$$ \partial _1a_{ij}=\partial _2a_{ij}=0,\tag 6 $$
$$\partial _1f_{ij}=f_{i-1j}a^{-1}_{ij} ,\tag 7$$
$$\partial _2f_{ij}=f_{i+1j}a_{i+1j} . $$

Consider Wronski matrices $W_k=W(f_{k1},\dots ,f_{kN})$,
using $\partial _2$ as a derivation and diagonal matrices 
$A_k={\text diag}(a_{k1},\dots , a_{kN})$ for $k=0,1,2,\dots ,n-1 $. 
Suppose that all these matrices are invertible and the quasideterminants
$|\partial _2W_k|_{NN}$ are defined. Set
$g_k=|\partial _2W_k|_{NN}|W_k|^{-1}_{1N}$, $k=0,1,\dots ,n-1$.

\proclaim{Theorem 5} Elements $g_k$, $k=0,1,\dots , n-1$
satisfy Toda system (1).
\endproclaim
{\bf Proof}. Set 
$\gamma _k= (\partial _2W_k)W_k^{-1}, k=0,1,\dots ,n-1$. Then 
apply Proposition 2, Lemma 3, and Proposition 1.
$\square $
\subhead 6. Toda systems over formal series of two variables
\endsubhead

To construct sets of elements $f_{ij}$, $a_{ij}$ satisfying 
the conditions of subsection 5 consider the algebra $Q=S[[u,v]]$ 
of formal series of commutative variables $u,v$ over a not
necessarily commutative algebra $S$.

Set $\partial _1={\partial \over {\partial u}}$,
    $\partial _2={\partial \over {\partial v}}$.

Choose any elements $a_{ij}\in S$.
Let us construct $f_{ij}$. Consider the diagonal matrices 
$a_j={\text diag}(a_{0j},\dots , a_{n-1j})$,
$j=1,\dots , N$. Let $R$ be the matrix of permutation
$\sigma (1,\dots,n)=(2,\dots, n, 1)$. Set $R_j=Ra_j$,
$j=1,\dots , N$.

Consider exponentials 
$e_j=e^{R_j^{-1}u + R_jv}$ as formal $n\times n$-matrix series of $u$
and $v$, $j=1,\dots , N$. For $j=1,\dots , N$ take any $1\times
n$-matrix $p_j$ (a row vector) over $S$ and set 
$f_j=p_je_j, j=1,\dots , N$ . Denote elements of $f_j$ by
 $f_{0j},\dots ,f_{n-1j}$.  
\proclaim {Proposition 6} Suppose that derivatives $\partial _1$,
$\partial _2$ are defined as above. Then elements $f_{ij}, i=0,1.\dots ,n-1,
j=1,\dots , N$ satisfy equalities (6) and (7).
\endproclaim
{\bf Proof}. It is enough to note that 
${\partial \over {\partial u}}e_j=e_jR_j^{-1}$, and
${\partial \over {\partial v}}e_j=e_jR_j$ for $j=1,\dots ,N$.
$\square $

\remark{Remark} Thus, we have constructed a family of solutions 
of the periodic Toda equations which are parametrized by $2nN$
$S$-valued parameters ($nN$ elements of the $n$-dimensional vectors 
$p_1,...,p_N$, and $nN$ elements $a_{ij}$). 
\endremark

\subhead 7. Example: the sine-Gordon equation \endsubhead
In the commutative case the sine-Gordon equation is 
obtained from the
2-periodic Toda system. We consider here the 2-periodic nonabelian
Toda system. In the notations of subsections 5,6 set

$$
f_{ij}=p_je^{a^{-1}_ju+a_jv} +(-1)^iq_je^{-a^{-1}_ju -a_jv},
i\in {\Bbb Z}/2{\Bbb Z},\ j=1,\dots , N. \tag 8
$$

with $p_j, q_j, a_j\in S$.

Then $(f_{ij})$ satisfy equations (6), (7) with $a_{ij}=a_j$.

Set $N=1$ and omit the index $j=1$ in our notations. According to
Theorem 5
$$
\gamma _i=(p-(-1)^iqe^{-2a^{-1}u-2av})a
(p+(-1)^iqe^{-2a^{-1}u-2av})^{-1},\ i\in {\Bbb Z}/2{\Bbb Z}
$$
are solutions of the 2-periodic Toda system.

\remark{Remark}If $S$ has an involutive automophism $*$
then $\gamma _0=\gamma ^*_1$ for appropriate $a,p,q$. 
Then $\gamma =\gamma _0$ satisfies
the noncommutative sine-Gordon equation
$$
\partial /\partial u((\partial \gamma /\partial v)\gamma ^{-1})=
\gamma (\gamma^*)^{-1}-\gamma^* \gamma ^{-1}.
$$

If $S=\Bbb C$ and $*$ is the complex conjugation then 
substitution $\gamma =e^{if}$ leads us to the sine-Gordon
equation.
\endremark 
For $N=2$ formulas (8) imply
$$
\gamma _i=(f_{i2}a^2_2-f_{i1}a_1f^{-1}_{i-1,1}f_{i-1,2}a_2)
(f_{i2}-f_{i1}a^{-1}_1f^{-1}_{i-1,1}f_{02}a_2)^{-1}=
$$
$$
=\{(p_2+(-1)^iq_2\eta_2)a_2 -(p_1+(-1)^iq_1\eta _1)a_1
(p_1-(-1)^iq_1\eta_1)^{-1}(p_2-(-1)^iq_2\eta_2)\}\times
$$
$$
\times \{(p_2+(-1)^iq_2\eta_2)a^{-1}_2 -(p_1+(-1)^iq_1\eta _1)a^{-1}_1
(p_1-(-1)^iq_1\eta_1)^{-1}(p_2-(-1)^iq_2\eta_2)\}^{-1},
$$

where $\eta _j=e^{-2a^{-1}_ju-2a_jv}, j=1,2 $.

\subhead 8. Langmuir equations and another form of Marchenko Lemma
\endsubhead

Let $R$ be a ring with unit, $\partial : R\rightarrow R$,
its derivation, and $\alpha :R\rightarrow R$ an
isomorphism, such that $\partial \alpha =\alpha \partial $.

Let also $\Gamma , A\in R$ be invertible elements such that
$\partial A=0$, $\partial \Gamma =\alpha ^2 (\Gamma )$,
and $\partial \Gamma +\Gamma = \alpha (\Gamma)A$.

Set $\gamma = (\partial \Gamma)\Gamma ^{-1}$. 
Since $\alpha $ is an isomorphism, $\gamma $ is invertible. 
The following lemma was proved in [M].

\proclaim {Lemma 7}The element 
$U=\gamma \alpha ^{-1} (\gamma) ^{-1}=1+\alpha (\gamma )-
\gamma $ satisfies the following
equality:
$$
\partial U=\alpha (U)U -U\alpha ^{-1} (U). \tag 9
$$
\endproclaim
 {\bf Proof}. We will show first that
$$
(\alpha (\gamma )-\gamma)(\gamma +1)=\partial \gamma , \tag 10
$$
$$
\alpha (\alpha (\gamma)-\gamma )\gamma = 
\alpha (\gamma )-\gamma . \tag 11
$$

Note, that $\partial \gamma = (\partial ^2\Gamma )\Gamma ^{-1}-
\gamma ^2$. From here one has
$$
\alpha (\gamma )=\alpha (\partial \Gamma)\alpha (\Gamma)^{-1}=
(\partial ^2\Gamma +\partial \Gamma)(\partial \Gamma +\Gamma)^{-1}=
$$
$$
=((\partial ^2\Gamma)\Gamma ^{-1}+\gamma)(\gamma +1)^{-1}=
(\partial \gamma +\gamma ^2 +\gamma)(\gamma +1)^{-1}=
\partial \gamma (\gamma +1)^{-1}+\gamma.
$$

This leads to (10). To prove (11), we notice that
$$
\alpha ^2(\gamma )=\alpha ^2(\partial \Gamma )\alpha ^2(\Gamma )^{-1}=
\partial (\alpha ^2(\Gamma ))\alpha ^2(\Gamma )^{-1}=
\partial ^2(\Gamma )(\partial \Gamma )^{-1}=
$$
$$
=(\partial ^2\Gamma)\Gamma ^{-1}\gamma ^{-1}=
(\partial \gamma +\gamma ^2)\gamma ^{-1}=
(\partial \gamma)\gamma ^{-1} +\gamma .
$$

Then 
$
(\alpha ^2(\gamma )-\alpha (\gamma ))\gamma =\partial \gamma
+\gamma ^2 -\alpha (\gamma )\gamma 
$.
By (10) the last expression equals to $\alpha (\gamma )-\gamma $
and (11) is proved.

Note that $\alpha (U)\gamma =\alpha (\gamma)$ and (10) implies

$$U=1+\alpha (\gamma ) -\gamma .\tag 12$$ 

Together with (10) it gives us
$U\alpha ^{-1} (U)=1+\partial \gamma $.

>From (10) it follows that $\partial (\alpha \gamma)=
\alpha ^3(\gamma)-\alpha (\gamma )$. Together with (11) it implies
$\alpha (U)U=1+\partial (\alpha \gamma )$.

>From (12) it follows that $\partial U=
\partial (\alpha \gamma )-\partial (\gamma )=
\alpha (U)U-U\alpha ^{-1} (U)$. $\square $

Similarly to subsection 3 we will consider the following
specialization of Lemma 7. 

Let $(Q,\partial) $ be a differential algebra 
with unit. Suppose we are given two sets of elements
$f_{ij}, a_{ij}$, $i\in \Bbb Z$, $j=1,2,\dots ,N$ in $Q$
such that for all $i,j$ $\partial a_{ij}=0$ and
$$
\partial f_{ij}=f_{i+2j},\ 
\partial f_{ij}+f_{ij}=f_{i+1j}a_{ij}.
\tag 13
$$

Suppose also that the Wronski matrices $W_k=W(f_{k1},\dots ,f_{kN})$
are invertible and the quasideterminants $|\partial W_k|_{NN}$
are defined and invertible.

Set 
$$
g_k=|\partial W_k|_{NN}|W_k|^{-1}_{1N}\cdot
|W_{k-1}|_{1N}|\partial W_{k-1}|^{-1}_{NN}=
$$
$$
=1+|\partial W_{k+1}|_{NN}|W_{k+1}|^{-1}_{1N}-
|\partial W_k|_{NN}|W_k|^{-1}_{1N}.
$$
\proclaim {Theorem 8} Elements $g_k, k=1,2,\dots $ satisfy generalized
Langmuir equations
$$
\partial g_k=g_{k+1}g_k-g_kg_{k-1}.
$$
\endproclaim

{\bf Proof}. The proof is similar to the proof of Theorem 5.
Set $R=Q^{\infty }=\{(w_k), k\in \Bbb Z\}$
be an algebra with coordinate-wise multiplication and 
coordinate-wise derivation induced by $\partial $. Set
$\alpha ((w_k))=(w_{k+1})$. Take $w_k=W_k$ and $A=(A_k)$,
$A_k={\text diag}(a_{k1},\dots ,a_{kN})$. Then apply
Lemma 7 and Proposition 1. $\square $

Note that in a commutative case the last equation can be written
in the form 
$$
\partial ({\text log}g_k)=g_{k+1}-g_{k-1}.
$$

{\bf Example}. Let $Q=S[[t]]$, where $S$ is an algebra with
unit and $t$ and independent variable commuting with elements
of $S$. Set $\partial =\partial /{\partial t}$. For 
$i=1,2,\dots , j=1,\dots , N$
consider a set of formal series
$$
f_{ij}(t)=p_j\mu ^i_je^{\mu ^2_jt}+q_j\mu ^{-i}_je^{\mu ^{-2}_jt},
$$

with $p_j, q_j, \mu _j \in S$.

These series satisfy equations (13) with $a_{ij}=\mu _j +\mu^{-1}_j$.
Let $N=1$. Omitting the index $j=1$ we get that functions
$$
g_k=(q+p\mu ^{2k+4}e^{(\mu ^2-\mu ^{-2})t})\mu ^{-2} 
(q+p\mu ^{2k}e^{(\mu ^2-\mu ^{-2})t})^{-1}\times
$$
$$
\times
(q+p\mu ^{2k-2}e^{(\mu ^2-\mu ^{-2})t})\mu ^2
(q+p\mu ^{2k+2}e^{(\mu ^2-\mu ^{-2})t})^{-1}
$$

satisfy the generalized Langmuir equations.

\remark{Remark} To obtain a solution of the periodic
Langmuir equations one has to take $\mu _j$'s to be the
corresponding roots of $1$.
\endremark
\subhead {9. Nonlinear Schr\"odinger Equation}\endsubhead

We will use the following form of the Marchenko Lemma. Let $R$
be an algebra with unit and let
$\partial ,\partial _0, \partial _1 $ be its derivations
commuting with each other. Let elements $A, B, C$ belong to the
common kernel of all three derivations and consider an
invertible element $\Gamma \in R$ such that
$$
\partial _0\Gamma +B\partial ^2\Gamma =\Gamma C,\tag 14
$$
$$
\partial _1\Gamma +B\partial \Gamma =\Gamma A.\tag 15
$$

The following lemma can also be found in [M].

\proclaim {Lemma 9} Let $\gamma =(\partial \Gamma)\Gamma ^{-1}$,
$U=\gamma B-B\gamma, V=\gamma B+B\gamma $. If $B^2=1$ then $U$ and $V$
satisfy the following equations
$$
2B\partial _0U +(\partial^2_1 +\partial ^2)U +2U^3
-2(U\partial _1V +(\partial _1V)U)=0, \tag 16
$$
$$
\partial _1V+B\partial V=U^2.\tag 17
$$
\endproclaim

{\bf Proof}. Multiply (14) by $\Gamma ^{-1}$ from the left, 
apply the derivation $\partial $, and then multiply by $\Gamma $
from the left. Then we get
$$
-\gamma \partial _0\Gamma + \partial _0\partial \Gamma -
\gamma B\partial ^2\Gamma +B\partial ^3\Gamma =0.\tag 18
$$
One can easily check that
$$
\partial _0\gamma=(\partial _0\partial \Gamma )\Gamma ^{-1}-
\gamma (\partial _0\Gamma )\Gamma ^{-1},
$$
$$
(\partial ^2\Gamma)\Gamma ^{-1}=\partial \gamma + \gamma ^2,
$$
$$
(\partial ^3\Gamma)\Gamma ^{-1}=\partial ^2\gamma 
+2(\partial \gamma )\gamma +\gamma \partial \gamma +\gamma ^3.
$$
The substitution of these formulas to (18) multiplied by 
$\Gamma ^{-1}$ from the right leads us to the equality
$$
\partial _0\gamma +B(\partial ^2\gamma +2(\partial \gamma)\gamma )
+[B,\gamma ](\partial \gamma +\gamma ^2)=0. \tag 19
$$
In a similar way (17) leads us to
$$
\partial _1 \gamma +B\partial \gamma +[B,\gamma]\gamma =0.\tag 20
$$  
Let us apply to (20) the derivation $\partial _1$. Then
$$
\partial ^2_1\gamma +B\partial \partial _1 \gamma +
[B,\partial _1 \gamma]\gamma +[B,\gamma ]\partial _1\gamma =0.
$$
Replacing in this formula $\partial _1\gamma $ by
$-B\partial \gamma -[B,\gamma ]\gamma $ and using the identity
$B^2=1$ one has
$$
(\partial ^2_1-\partial ^2)\gamma - 2(B[B,\partial \gamma ]\gamma
+ B[B,\gamma ]\gamma ^2 +[B,\gamma ]^2\gamma )=0.
$$
It allows us to exclude from (19) the term $[B,\gamma ]\gamma ^2$
and get
$$
2\partial _0\gamma +B(\partial ^2 +\partial ^2_1)\gamma
+4B(\partial \gamma)\gamma +2[B,\gamma ]\partial \gamma
-2[B,\partial \gamma ]\gamma -2B[B,\gamma ]^2\gamma =0,
$$
or
$$
2\partial _0\gamma +B(\partial ^2 +\partial ^2_1)\gamma
+2\{B,\partial \gamma \}\gamma 
+2[B,\gamma ](\partial \gamma -B[B,\gamma ]\gamma )=0,\tag 21
$$
where $\{x,y\}=xy+yx$.

Multiplication of (21) by $B$ from the left and using (20)
leads us to
$$
2B\partial _0 \gamma +(\partial ^2+\partial ^2_1)\gamma
+2[B,\gamma ]^2\gamma -2\{B,\partial _1\gamma \}\gamma
+2[B,\gamma ]\partial _1\gamma =0.
$$

Commutation the last equality with $B$ gives us (16). Note also,
that (20) implies (17). $\square $

We will apply Lemma 9 when $\partial _1=0$. It will give us a
nonabelian version of the nonlinear Schr\"odinger equation. 
When $\partial _1\neq 0$ similar arguments lead to the nonabelian 
Davey-Stewartson system (cf.[LY]). 

Let $Q$ be an algebra with unit and with two commuting derivations
$\partial _0$ and $\partial $. Let $R$ be the algebra of $N\times N$-
matrices over $Q$. We extend in an obvious way the derivations 
$\partial _0$ and $\partial $
to the derivations of $R$ and denote them by the same letters.
For an element $f\in Q$
denote $\partial ^kf$ by $f^{(k)}$, and for a set $f_1,\dots , f_N$
let $W=W(f_1,\dots ,f_N)=(f^{(k-1)}_j),
k,j=1,\dots , N$. Suppose that $W$ is invertible and that 
$N\times N$-matrices $\Gamma =W, B={\text diag}
\{b,\dots , b\}, A$ satisfy equations (14) in $R$. Suppose also
that $b^2=1$ and
set $g=|\partial W|_{NN}|W|^{-1}_{1N}$. Then Lemma 9 implies

\proclaim {Theorem 10} The element $U=gb-bg$ satisfies the equation

$$2b\partial _0U+\partial ^2U +2U^3=0.\tag 22$$

\endproclaim

If $Q$ itself is an algebra of $r\times r$-matrices over an algebra
with unit and $b={\text diag}\{1,\dots ,1, -1,\dots , -1\}$
(there are $r_1$ 1's and $r_2=r-r_1$ (-1)'s) then $U=(U_{pq})$
$p,q=1,2$ is a block-matrix where $U_{pq}$ is an
$r_p\times r_q$-matrix
and $U_{11}=0, U_{22}=0$. Then (22) implies

$$
2\partial _0U_{12}+\partial ^2 U_{12} +2U_{12}U_{21}U_{12}=0,\tag 23
$$
$$
-2\partial _0U_{21}+\partial ^2 U_{21} +2U_{21}U_{12}U_{21}=0.\tag 24
$$

Let $Q$ is the algebra of matrices over an algebra with an
involutive automorphism $*$. Denote by $X^*$ the hermitian
conjugation to the matrix $X$. If
$r_1=r_2$ then one could choose $A$ such that
$U_{21}=U_{12}^*$. From $*\partial _0*=-\partial _0 $,
$*\partial *=\partial $ and (23) we get an equation for $X=U_{12}$:
$$
2\partial _0X+\partial ^2X +2XX^*X=0.
$$ 

As an example we will consider an algebra of formal series
$R[[u,v]]$ over a not necessarily commutative algebra $R$ with unit over
complex numbers. Here variables $u$ and $v$ commute with each other
and with elements of $R$. We set $\partial _0=i\partial /\partial u,
\partial =\partial /\partial v$, where $i=\sqrt {-1}$. Let $b^2=1, b\in R$,
set $q_1=(1+b)/2, q_2=(1-b)/2$. Then $bq_1=q_1, bq_2=-q_2$.
Let $a_1,\dots ,a_N,
c_1,\dots ,c_N, d_1, \dots , d_N\in R$. Set
$$
f_j=q_1c_je^{a_jv+ia^2_ju} +q_2d_je^{-a_jv-ia^2_ju}, 
\ j=1,\dots ,N.
$$
It is easy to see that $W=W(f_1,\dots ,f_N)$ satisfies
equations (14) with $C=0, A={\text diag} \{a_1,\dots, a_N\}$.

As an example, we consider case $N=1$. Then $U=gb-bg$ where
$$
g= (q_1ce^{av+ia^2u} +q_2de^{-av -ia^2u})a
(q_1ce^{au+ia^2v} +q_2de^{-av-ia^2u})^{-1}.
$$

is a solution of equation (22).

One may consider a case when $R$ is the algebra of 
$2\times 2$-matrices over a $\Bbb C$-algebra $R_0$ with an
involutive automorphism $a\rightarrow a^*$ which agrees with the complex
conjugation. For appropriate $c$ and $d$
$f\in R[[u,v]]$ could be written as
$$
f=\left(\matrix \alpha e^x & \beta ^* e^{-x^*}\\
        \beta e^{-x} & \alpha ^* e^{x^*}\endmatrix\right),
$$

where $x=av+ia^2u$, and  $a,\alpha ,\beta \in R_0$.

Then $\Gamma =f$ satisfies equations (14) with $C=0$,
$B={\text diag}\{1,-1\}, A={\text diag}\{a, -a^*\}$.

Consider the $2\times 2$-matrix 
$(\gamma _{pq})=(\partial f/\partial v)f^{-1}$.
One can see that $\gamma _{21}=\gamma ^* _{12}$. It follows
that 
$$
w=-2\gamma _{12}= 2(\alpha a\alpha ^{-1}+
\beta ^*a^*(\beta ^*)^{-1})
(\alpha ^*e^{2x^*}(\beta ^*)^{-1} -
\beta e^{2x}\alpha ^{-1})^{-1}  
$$
satisfies the equation

$$
2\partial _0w+\partial ^2w+2ww^*w=0.
$$

When $R_0=\Bbb C$ one could write 
$\alpha e^{2x}\beta ^{-1}=e^y$ and then  
arrive to a well-known example [LS]
$$
w=(a+\bar a)e^{-iI}/{{\text sinh}\ R},
$$
where $I={\text Im}\ y, R={\text Re}\ y$.
\remark{Remark} One could change the variable $U\rightarrow iU$
in the Schr\"odinger equation and arrive to the nonlinear heat equation.
Then we could eliminate $i=\sqrt {-1}$ in our formulas and to generalize them
to the case of any field of characteristic zero.
\endremark


\Refs
\ref\by [EGR] Etingof, P., Gelfand, I., and Retakh, V.\paper
Factorization of differential operators, quasideterminants,
and nonabelian Toda field equations \jour Math. Research Letters
\vol 4\yr 1997\pages 413-425
\endref

\ref\by [GR1] Gelfand, I., and Retakh, V.\paper 
Determinants of matrices over noncommutative rings
\jour Funct.An. Appl.\vol 25\issue 2\yr 1991\pages 91-102
\endref

\ref\by [GR2] Gelfand, I., and Retakh, V.\paper 
A theory of noncommutative determinants and characteristic functions 
of graphs
\jour Funct.An. Appl.\vol 26\issue 4\pages 1-20\yr 1992
\endref

\ref\by [GR3] Gelfand, I., and Retakh, V.\paper
A theory of noncommutative determinants and characteristic functions
of graphs.I
\jour in: Publ. LACIM, UQAM\vol 14\yr 1993\pages 1-26\endref

\ref\by [GR4] Gelfand, I., and Retakh V.\paper Quasideterminants, I
\jour Selecta Math. (to appear) \endref

\ref\by [Kr] Krichever, I.M.\paper The periodic non-abelian Toda chain 
and its two-dimensional generalization\jour Russ. Math. Surv.\vol 36
\issue 2\yr 1981\pages 82-89\endref

\ref\by [Kr] Krichever, I.M., Novikov, S.P.
\paper Holomorphic bundles over algebraic curves 
and nonlinear equations\jour Russ. Math. Surv.\vol 35
\issue 6\yr 1980\pages 53-79\endref

\ref\by [LS] Leznov, A.N., Saveliev M.V.\book 
Group Methods of Integration of Nonlinear Dynamical Systems
(in Russian)\publ Nauka\publaddr Moscow\yr 1985\endref

\ref\by [LY] Leznov, A.N., Yuzbashyan, E.A.\paper Multi-soliton
solutions of Two-dimensional Davey-Stewartson Equation
\jour hep-th 9612107\endref

\ref\by [M] Marchenko V.A. \book Nonlinear equations and operator
algebras (in Russian)\publ Naukova Dumka \publaddr Kiev \yr 1986
\endref

\ref\by [P] Perk, J.H.H.\paper Equations of motions for the transverse
correlations of the one-dimensional $XYZ$-model at finite temperature
\jour Phys. Lett. A\vol 92\yr 1980 \pages 1-5\endref

\ref\by [PC] Perk, J.H.H. and Capel, H.W.\paper Transverse correlations
in the inhomogeneous $XY$-model at infinite temperature \jour Physica A
\vol 92 \yr 1978 \pages 163-184\endref

\ref\by [RS] Razumov, A.V., Saveliev, M.V.\paper Maximally Nonabelian Toda
systems\jour Nuclear Physics B, hep-th 9612081\endref

\endRefs
\enddocument